\documentclass[11pt,reqno,tbtags]{amsart}
\usepackage{amssymb}
\usepackage{natbib}
\bibpunct[, ]{[}{]}{;}{n}{,}{,}

\numberwithin{equation}{section}
\allowdisplaybreaks

\newtheorem*{property*}{Property \csname @currentlabel\endcsname}

\newtheorem{theorem}{Theorem}[section]
\newtheorem{lemma}[theorem]{Lemma}

\newtheorem{corollary}[theorem]{Corollary}

\theoremstyle{definition}
\newtheorem{example}[theorem]{Example}
\newtheorem{definition}[theorem]{Definition}

\newtheorem{remark}[theorem]{Remark}

\theoremstyle{remark}

\newenvironment{romenumerate}{\begin{enumerate}
 }{\end{enumerate}}

\newcounter{oldenumi}
{\setcounter{oldenumi}{\value{enumi}}
\begin{romenumerate} \setcounter{enumi}{\value{oldenumi}}}
{\end{romenumerate}}

\newcounter{thmenumerate}

\newcounter{xenumerate}   

\newcommand\pfitem[1]{\par(#1):}

\newcommand{\refT}[1]{Theorem~\ref{#1}}
\newcommand{\refC}[1]{Corollary~\ref{#1}}
\newcommand{\refL}[1]{Lemma~\ref{#1}}
\newcommand{\refR}[1]{Remark~\ref{#1}}
\newcommand{\refS}[1]{Section~\ref{#1}}
\newcommand{\refD}[1]{Definition~\ref{#1}}

\newcommand{\refE}[1]{Example~\ref{#1}}

\newcommand{\refand}[2]{\ref{#1} and~\ref{#2}}

\newcommand\marginal[1]{\marginpar{\raggedright\parindent=0pt\tiny #1}}

\begingroup
  \divide\count255 by 60
  \multiply\count255 by -60
  \advance\count255 by \time
  \ifnum \count255 < 10 \xdef\klockan{\the\count1.0\the\count255}
  \else\xdef\klockan{\the\count1.\the\count255}\fi
\endgroup

\newcommand{\sumin}{\sum_{i=1}^n}
\newcommand{\sumiN}{\sum_{i=1}^N}

\newcommand{\prodiN}{\prod_{i=1}^N}
\newcommand{\prodoo}{\prod_{i=1}^\infty}

\newcommand\set[1]{\ensuremath{\{#1\}}}

\newcommand\bigpar[1]{\bigl(#1\bigr)}
\newcommand\Bigpar[1]{\Bigl(#1\Bigr)}
\newcommand\biggpar[1]{\biggl(#1\biggr)}

\newcommand\lrpar[1]{\left(#1\right)}

\newcommand\lrabs[1]{\left|#1\right|}
\def\rompar(#1){\textup(#1\textup)}    

\newcommand\parfrac[2]{\Bigpar{\frac{#1}{#2}}}

\def\xexp(#1){e^{#1}}

\newcommand\ntoo{\ensuremath{{n\to\infty}}}

\newcommand\ttoo{\ensuremath{{t\to\infty}}}
\newcommand\asttoo{\text{as }\ttoo}
\newcommand\bmin{\wedge}
\newcommand\bmax{\vee}
\newcommand\norm[1]{\|#1\|}

\newcommand\iid{i.i.d.\spacefactor=1000}    
\newcommand\ie{i.e.\spacefactor=1000}
\newcommand\eg{e.g.\spacefactor=1000}

\newcommand\cf{cf.\spacefactor=1000}
\newcommand{\as}{a.s.\spacefactor=1000}

\newcommand\whp{{whp}}

\newcommand{\tend}{\longrightarrow}
\newcommand\dto{\overset{\mathrm{d}}{\tend}}
\newcommand\pto{\overset{\mathrm{p}}{\tend}}

\newcommand\eqd{\overset{\mathrm{d}}{=}}

\newcommand\bbR{\mathbb R}

\newcommand\bbN{\mathbb N}

\newcommand\bbZ{\mathbb Z}

\newcounter{CC}
\newcounter{cc}

\newcommand\E{\operatorname{\mathbb E{}}}
\renewcommand\P{\operatorname{\mathbb P{}}}

\newcommand\Po{\operatorname{Po}}
\newcommand\Bi{\operatorname{Bi}}
\newcommand\Be{\operatorname{Be}}

\newcommand\N{\operatorname{N}}

\newcommand\ga{\alpha}

\newcommand\gd{\delta}

\newcommand\gf{\varphi}

\newcommand\gk{\kappa}
\newcommand\gl{\lambda}
\newcommand\gL{\Lambda}

\newcommand\gs{\sigma}

\newcommand\eps{\varepsilon}

\newcommand\cA{\mathcal A}

\newcommand\cE{\mathcal E}
\newcommand\cF{\mathcal F}
\newcommand\cL{{\mathcal L}}
\newcommand\cS{{\mathcal S}}

\newcommand\cX{{\mathcal X}}

\newcommand\ett[1]{\boldsymbol1[#1]} 

\def\[#1]{[\![#1]\!]}

\newcommand\qq{^{1/2}}
\newcommand\qqw{^{-1/2}}
\newcommand\qw{^{-1}}
\newcommand\qww{^{-2}}

\renewcommand{\=}{:=}

\newcommand\oi{[0,1]}
\newcommand\oiset{\set{0,1}}
\newcommand\ooo{[0,\infty)}

\newcommand\dtv{d_{\mathrm{TV}}}
\renewcommand\dh{d_{\mathrm{H}}}
\newcommand\dhh{d_{\mathrm{H}}^2}

\newcommand\lhs{left-hand side}

\newcommand\aseq{\cong}
\newcommand\aseqq{asymptotically equivalent}
\newcommand\cont{\vartriangleleft}
\newcommand\contr{\vartriangleright}
\newcommand\contxx{\cont\contr}
\newcommand\iij{I_{ij}}
\newcommand\ini{I_{ni}}
\newcommand\pij{p_{ij}}
\newcommand\pijn{p_{ij}}
\newcommand\pijny{\pijn'}
\newcommand\pijx{\set{\pij}}

\newcommand\pijnx{\set{\pijn}}
\newcommand\pijnxy{\set{\pijn'}}
\newcommand\pni{p_{ni}}
\newcommand\pnix{\set{\pni}}
\newcommand\pnixy{\set{\pni'}}
\newcommand\qni{q_{ni}}
\newcommand\qij{q_{ij}}
\newcommand\Pni{P_{ni}}
\newcommand\Pn{P_{n}}
\newcommand\gnpij{\ensuremath{G(n,\pijx)}}

\newcommand\gnpijn{\ensuremath{G(n,\pijnx)}}
\newcommand\gnpijny{\ensuremath{G(n,\pijnxy)}}
\newcommand\gnp{\ensuremath{G(n,p)}}
\newcommand\gnpy{\ensuremath{G(n,p')}}
\newcommand\pp{{\mathbf p}}
\newcommand\gnpp{\ensuremath{G(n,\pp)}}

\newcommand\gnppy{\ensuremath{G(n,\pp')}}
\newcommand\nn{[n]}
\newcommand\ijn{{1\le i<j\le n}}
\newcommand\equ{\asymp}
\newcommand\xiin{_{i=1}^{N}}
\newcommand\xiinn{_{i=1}^{N(n)}}
\newcommand\suminn{\sum_{i=1}^{N(n)}}
\newcommand\xppn{(X_n\mid\pp_n)}
\newcommand\xppny{(X_n'\mid\pp_n')}
\newcommand\xppnyy{(X_n''\mid\pp_n'')}
\newcommand\xppnyyy{(X_n'''\mid\pp_n''')}
\newcommand\bp{\overline P}
\newcommand\bq{\overline Q}
\newcommand\prodn{\prod_1^n}
\newcommand\prodNn{\prod_1^{N(n)}}
\newcommand\sumNn{\sum_1^{N(n)}}
\newcommand\Bex{\Be}
\newcommand\supx[1]{^{(#1)}}
\newcommand\pijq[1]{\pij\supx{#1}}

\newcommand\pijo{\hat p_{ij}}
\newcommand\cSS{\cS}
\newcommand\cXX{\cX}
\newcommand\gdn{\gd_n}

\newcommand\CS{Cauchy--Schwarz}
\newcommand\CSineq{\CS{} inequality}

\newcommand\REM[1]{{\raggedright\texttt{[#1]}\par\marginal{XXX}}}

\newcommand\urladdrx[1]{{\urladdr{\def~{{\tiny$\sim$}}#1}}}

\begin{document}
\title[Equivalence and contiguity of random graphs]
{Asymptotic equivalence and contiguity of some random graphs}

\date{February 12, 2008} 

\author{Svante Janson}
\address{Department of Mathematics, Uppsala University, PO Box 480,
SE-751~06 Uppsala, Sweden}
\email{svante.janson@math.uu.se}
\urladdrx{http://www.math.uu.se/~svante/}

\subjclass[2000]{} 

\begin{abstract} 
We show that asymptotic equivalence, in a strong form, holds between
two random graph models with slightly differing edge probabilities
under substantially weaker conditions than what might naively be
expected.

One application is a simple proof of a recent result by
van den Esker, van der Hofstad and Hooghiemstra     
on the equivalence between graph distances for some random graph models.
\end{abstract}

\maketitle

\section{Introduction}\label{S:intro}

There are many different models of random graphs. Sometimes, the
differences are minor, and it can be guessed that the asymptotic
behaviour of two models are the same (for all or at least for some
interesting properties). This note concerns some cases where it is
possible to actually prove such results in a strong form. We begin by
defining the two types of asymptotic equality that we will study.
All unspecified limits are as \ntoo.

\begin{definition}\label{D1}
  Let $(\cXX_n,\cA_n)$, $n\ge1$, be a sequence of arbitrary measurable
  spaces
and let $P_n$ and $Q_n$ be two probability measures on
  $(\cXX_n,\cA_n)$.
\begin{romenumerate}
  \item
The sequence $(P_n)_n$ is \emph{\aseqq}
to
$(Q_n)_n$, denoted by $(P_n)_n\aseq(Q_n)_n$, 
if for every sequence of measurable sets $A_n$ (\ie,
$A_n\in\cA_n$), we have $P_n(A_n)-Q_n(A_n)\to0$.
  \item
The sequence $(P_n)_n$ is \emph{contiguous} with respect to
$(Q_n)_n$,
denoted by $(P_n)_n\cont(Q_n)_n$, 
if for every sequence of measurable sets $A_n$ 
such that $Q_n(A_n)\to0$, we also have $P_n(A_n)\to0$.
\end{romenumerate}
We use the same terminology and notations for sequences of random variables
$X_n$ and $Y_n$ with values in the same space $\cXX_n$, meaning
that these properties hold for their 
distributions $\cL(X_n)$ and $\cL(Y_n)$.  
For example, $(X_n)_n\aseq(Y_n)_n$ means
that
$\P(X_n\in A_n)-\P(Y_n\in A_n)\to0$ for every sequence $(A_n)_n$.
We will also use the simpler notations $X_n\aseq Y_n$ and $X_n\cont
Y_n$, etc.
\end{definition}

Note that asymptotic equivalence is a symmetric relation while contiguity
is not; we say that $(P_n)_n$ and $(Q_n)_n$ 
are (mutually) contiguous, $(P_n)_n\contxx(Q_n)_n$, if 
both $(P_n)_n\cont(Q_n)_n$ and $(Q_n)_n\cont(P_n)_n$, \ie,
if $P_n(A_n)\to0\iff Q_n(A_n)\to0$ for any sequence of measurable sets
$A_n\subseteq\cXX_n$. 
(And similarly for sequences of random variables  $X_n$ and $Y_n$.)

Asymptotic equivalence implies contiguity, but not conversely
(see \eg{} \refE{EA1} and \refR{Rreg}), so
contiguity is a weaker property.

We illustrate 
these notions by two simple examples.

\begin{example}\label{EA1}
In the special case of
two constant sequences, $P_n=P$ and $Q_n=Q$ where $P$ and $Q$ are two
probability measures defined on the same space  $(\cXX_n,\cA_n)=(\cXX;\cA)$, 
$(P_n)_n\aseq(Q_n)_n$ if and only
  $P=Q$, and 
$(P_n)_n\cont(Q_n)_n$ if and only if $P\ll Q$, \ie, $P$
  is absolutely continuous with respect to $Q$. Hence asymptotic
  equivalence and contiguity
  can be thought of as asymptotic versions of equality
and absolute continuity, respectively.
\end{example}

\begin{example}\label{EA2}
Let $X_n$ be random elements in some spaces $\cXX_n$ and let $\cE_n$ be
events that depend on $X_n$ only, \ie, $\cE_n=\set{X_n\in E_n}$ for
some (measurable) sets $E_n\subseteq\cXX_n$, and suppose that
  $\liminf\P(\cE_n)>0$.
Let $Y_n\=(X_n\mid \cE_n)$ be $X_n$ conditioned on $\cE_n$ (possibly
ignoring some small $n$ with $\P(\cE_n)=0$).
Then, for any $A_n$, and some $C<\infty$,
\begin{equation*}
 \P(Y_n\in A_n)=\frac{\P(X_n\in A_n\cap E_n)}{\P(\cE_n)}
\le C \P(X_n\in A_n\cap E_n)
\le C \P(X_n\in A_n)
\end{equation*}
and thus $(Y_n)_n\cont (X_n)_n$

An important random graph example of this is when $Y_n$ is a random
graph with a given degree sequence $d_1,\dots,d_n$, uniformly chosen
among all such graphs, and $X_n$ is the random multigraph constructed
by the configuration model (see, \eg, \citet{Bollobas}); then
$Y_n\eqd(X_n\mid X_n\text{ is a simple graph})$ so 
$(Y_n)_n\cont(X_n)_n$ provided $\liminf_\ntoo\P(X_n\text{ is
simple})>0$.
This is the case when $\sumin d_i\to\infty$ and $\sumin
d_i^2=O\bigpar{\sumin d_i}$, see \citet{SJ195} (with several earlier partial
results by various authors), which makes it possible to transfer many
results from $X_n$ to $Y_n$. Indeed, this is a standard method to
study random graphs with a given degree sequence, and in particular
random regular graphs, see \eg, \cite{Bollobas}, \cite{JLR}, \cite{Wormald}.
\end{example}


\begin{remark}\label{R1}
  Suppose that $X_n\aseq Y_n$. If $\P(X_n\in A_n)\to\ga$ for some
  sequence of (measurable) sets $A_n$ and some $\ga\in\oi$, then 
$\P(Y_n\in A_n)\to\ga$ too. Hence, any result for $X_n$ that can be stated
  in terms of convergence of some probabilities holds for $Y_n$ too;
  for example, this includes any result of the type $\gf_n(X_n)\pto a$
  and $\gf_n(X_n)\dto W$ for some functionals $\gf_n:\cXX_n\to\bbR$ (and
  a number $a$ or a random variable $W$).
However, results that are sensitive to events with small
  probabilities, such as moment convergence or large deviation
  estimates, do not transfer automatically. For example, if 
$X_n\aseq Y_n$ and we
  know that $\E\gf_n(X_n)\to a$, we may guess that $\E\gf_n(Y_n)\to a$
  too, but we cannot conclude it without further information (for
  example uniform integrability of $\gf_n(X_n)$ and $\gf_n(Y_n)$). 

If instead only $X_n\contr Y_n$, then results of the type
$\gf_n(X_n)\pto a$ still transfer to $Y_n$, but not result on
convergence in distribution. (If $\gf_n(X_n)\dto W$, then the sequence
$\gf_n(Y_n)$ is tight, but does not have to converge to $W$, or at
all. Typically, $\gf(Y_n)\dto W'$ for some $W'\neq W$, see \eg{}
\refE{EA1} and several examples of cycle counts
in \cite[Chapter 9]{JLR}.)
\end{remark}

Suppose now that $G_n$ and $G'_n$ are random graphs on the vertex set
$[n]\=\set{1,\dots,n}$. 
By the standard \refT{T0} below, $G_n\aseq G'_n$ if and only if it is
possible to couple $G_n$ and $G'_n$, \ie, to define them
simultaneously on some probability space, such that 
$\P(G_n\neq G_n')\to0$. (We assume that we are interested only
in the distributions of $G_n$ and $G_n'$, so we may replace them by
any random graphs with the same distributions.)

In particular, we will study random graphs of the following type.
If $\pij$, $1\le i<j\le n$, are given probabilities in 
\oi, let $\gnpij$ be the random graph on $[n]$ where the edge $ij$
appears with probability $\pij$ and the
indicators $\iij\=\ett{\text{edge $ij$ appears}}$, $1\le i<j\le n$, 
are independent.
(We will later also consider an extension to random $\pij$, see
\refS{Sresults}.) 

Consider two sequences of such graphs, defined by probabilities
$\set{\pijn}_{1\le i<j\le n}$ and
$\set{\pijn'}_{1\le i<j\le n}$;  $\pij$ and $\pijn$ may depend
on $n$ too, but to simplify the notation we do not show this explicitly. 
It is obvious that we may couple the edge indicators $\iij$ and
$\iij'$ of $ij$ in
$\gnpijn$ and $\gnpijny$ such that
$\P(\iij\neq\iij')=|\pijn-\pijn'|$, and by taking independent pairs
$(\iij,\iij')$ we obtain a coupling of the random graphs
$\gnpijn$ and $\gnpijny$ with
\begin{equation}\label{simple}
 \P\bigpar{\gnpijn\neq\gnpijny}\le\sum_{i<j}|\pijn-\pijn'|. 
\end{equation}
Consequently, $\gnpijn\aseq\gnpijny$ if 
$\sum_{i<j}|\pijn-\pijn'|\to0$; a simple fact that has been used by
many authors. It may be believed that this is essentially best
possible, but, somewhat surprisingly, this is not so.
In fact, by \refC{C1} below, provided $\pijn\le0.9$, say,
$\gnpijn\aseq\gnpijny$ if 
$\sum_{i<j}(\pijn-\pijn')^2/\pijn\to0$. 
(Moreover, \refT{T1}(i) shows that this is best possible
if, for example, $\pijny\le2\pijn$.)

For a particular case,
suppose that $\pijn'=\pijn+O(\pijn^2)$, 
Then \refC{C2} shows that $\gnpijn\aseq\gnpijny$ if 
$\sum_{i<j}\pijn^3\to0$, while \eqref{simple} implies this only under
the stronger condition
$\sum_{i<j}\pijn^2\to0$.
A typical case where this is an important improvement is when all
$\pijn=\Theta(1/n)$ and $|\pijn'-\pijn|=\Theta(1/n^2)$.
See further the examples in \refS{Sex}.

To see that such an improvement of \eqref{simple} might be possible at
all, consider
as an example the case when all $\pijn$ are the same, so we consider
the random graph $\gnp$:

\begin{example}\label{Egnp}
Let $p=p(n)$ and $p'=p'(n)$ be given in $\oi$ and consider $\gnp$ and $\gnpy$.
Let $N\=\binom n2$ be the number of possible
edges and let $M \sim \Bi(N,p)$ and $M'\sim\Bi(N,p')$ be the number of
edges in $\gnp$ and $\gnpy$. The conditional distribution of $\gnp$
given $M=m$ is uniform over all graphs on $\nn$ with $m$ edges, and
the conditional distribution of $\gnpy$ given $M'=m$ is the same.
It follows that any coupling of $M$ and $M'$ may be extended to a
coupling of $\gnp$ and $\gnpy$ such that $\gnp=\gnpy$ when $M=M'$; as
a consequence, using \eqref{coupling} below,
$\dtv\bigpar{\gnp,\gnpy}=\dtv(M,M')$, and in particular, as \ntoo,
using also
\refT{T0}, 
$$\gnp\aseq\gnpy\iff M\aseq M'.$$

Since $M$ and $M'$ have binomial distributions with the same $n$,
$\P(M=k)/\P(M'=k)$ is monotone in $k$, and it follows that the maximum
of $|\P(M\in A)-\P(M'\in A)|$ over subsets $A$ of $\bbZ$ is attained for
a set of the form $[0,1,\dots,j]$. 

Now suppose that $p\to0$ and
$N(p'-p)/\sqrt{Np}\to\ga\in[-\infty,\infty]$.
Suppose first that $\ga$ is finite.
By the central limit theorem,
$(M-Np)/\sqrt{Np}\dto \N(0,1)$ and  
$(M'-Np)/\sqrt{Np}\dto \N(\ga,1)$, and it follows 
easily that
\begin{equation*}
  \begin{split}
 \dtv\bigpar{\gnp,\gnpy}
&=\dtv(M,M')
=\sup_j|\P(M\le j)-\P(M'\le j)|
\\&
=\sup_x|\P\bigpar{\N(0,1)\le x}-\P\bigpar{\N(\ga,1)\le x}|
+o(1)
\\&
\to \Phi(\ga/2)-\Phi(-\ga/2).	
  \end{split}
\end{equation*}
It follows, using \refT{T0} again,  
that $\gnp\aseq\gnpy$ if and only if $\ga=0$, \ie{}
$N(p'-p)/\sqrt{Np}\to0$, which is equivalent to 
$\sum_{i<j}(p'-p)^2/p=N(p'-p)^2/p\to0$. 
For example, if $p=1/n$ and $p'=1-e^{-1/n}=p-\frac12n\qww+O(n^{-3})$, 
then $N(p'-p)^2/p =O(1/n)$, so $\gnp\aseq\gnpy$, but
$N|p-p'|\to1/4$, so \eqref{simple} is not enough to show this.

We see that in this example, the trick to improve the simple and
'obvious' edgewise coupling used in \eqref{simple}, is to first ignore the
positions of the edges and couple their numbers only; this then is
extended to a coupling of the random graphs by randomly reinserting
the positions. \refC{C1}(i) shows that couplings improving the simple edgewise
coupling exist also when the edge probabilities are unequal, but in
that case we do not know any explicit construction of such couplings.
\end{example}

We give the main results in \refS{Sresults} and a number of
examples in \refS{Sex}; this includes an application to a recent
result by
van den Esker, van der Hofstad and Hooghiemstra   
(\refE{EEHH}).
Proofs are given in 
\refS{Spf}, after some preliminaries in \refS{Smore}.

We use the standard notations $o_p$ and $O_p$, 
see \eg{} \cite[Section 1.2]{JLR}, and 
we write {\whp} 
(with high probability) for events with probability tending to 1 as \ntoo.

\begin{remark}  \label{Rreg}
  There are also interesting examples of contiguity among random
  graphs of other types than \gnpij. In particular, several different
  constructions of random regular graphs (or multigraphs) are known to
  yield distributions that are (mutually) contiguous but not \aseqq,
  see \eg{} \cite{SJ102}, 
  \cite[Section 9.5]{JLR}, 
  \cite{SJ140}.
  These examples are not covered by the present paper.
\end{remark}

\section{Results} \label{Sresults}

We defined above the random graph \gnpp, where $\pp=\pijnx_\ijn$ is a
vector of probabilities. We extend the definition 
of \gnpp{}
to the case when
$\pp$ is a random vector (with entries in \oi) by conditioning on
$\pp$, 
\ie, given $\pp=\pijnx$, the edge indicators $\iij$ are independent
with $\iij\sim\Be(\pij)$. 
Random graphs of this type have been studied in many papers, see for example
\citet{SJ178} and the further references given there.

We now state our main results on \aseqq{} and contiguity of such
random graphs. Actually, the results have nothing to do with the graph
structure and the way the indicator variables are indexed by pairs $ij$. It
therefore seems more natural to consider the more general situation of
a (finite or infinite) sequence $(I_i)_1^N$ of indicator
variables. (An indicator variable is a random variable with values in
$\oiset$, \ie{} a random variable with a Bernoulli distribution
$\Be(p)$ for some $p\in\oi$.) The results for random graphs then
follow by relabelling the indicators.

We define a function $\rho:\oi^2\to\ooo$ in \refD{Drho}, where
we also give some equivalent (within constant factors) alternative
formulas that often are more convenient. Since the results below are
not affected by changing $\rho$ within constant factors, we could use
any of these alternative formulas 
(and several other similar ones) as our definition.
(The motivation for the definition comes in \refL{LBe}.)

We write $x\equ y$ (where $x,y\ge0$)
to denote that $cy\le x\le Cy$ for some positive constants $c,C$, \ie,
that $x=\Theta(y)$  (or, equivalently, $x=O(y)$ and
$y=O(x)$). Further, we use $x\bmax y$ and $x\bmin y$ for the maximum and
minimum, respectively, of $x$ and $y$.
We interpret $0/0$ as 0.

\begin{definition}
\label{Drho}
Let
  \begin{align}
\rho(p,q)
&\= \bigpar{\sqrt p-\sqrt{q}}^2 + \bigpar{\sqrt{1-p}-\sqrt{1-q}}^2
\label{rho1} \\
&\equ 
\frac{(p-q)^2}{p+q} + \frac{(p-q)^2}{1-p+1-q}
\label{rho2} \\
&\equ 
\frac{(p-q)^2}{(p\bmax q)\bmin((1-p)\bmax(1-q))}
\label{rho3} \\
&\equ 
\frac{(p-q)^2}{p\bmin(1-p)} \bmin|p-q|.
\label{rho4} 
\intertext{In particular, if $p\le0.9$, then}
\rho(p,q)
&\equ 
\frac{(p-q)^2}{p} \bmin|p-q|.
\label{rho11} 
  \end{align}  
\end{definition}
Of course, the constant $0.9$ here and below is arbitrary and could be
replaced by any number $<1$.

\begin{proof}
  The first equivalence follows from
  \begin{equation*}
\bigpar{\sqrt{p}-\sqrt{q}}^2 
= \frac{(p-q)^2}{(\sqrt{p}+\sqrt{q})^2}	
\equ \frac{(p-q)^2}{p+q},
  \end{equation*}
together with the similar result with $1-p$ and $1-q$.
The second follows from $x+y\equ x\bmax y$ for $x,y\ge0$ (used
thrice).
The third is equivalent to
\begin{equation}\label{lrho}
  (p\bmax q)\bmin\bigpar{(1-p)\bmax(1-q)}
\equ
\bigpar{p\bmin(1-p)}\bmax|p-q|,
\end{equation}
which is easily verified, for example by assuming (by the symmetry
$p\mapsto1-p$, $q\mapsto1-q$) that $p\le1/2$, in which case
\eqref{lrho} easily reduces to 
$p\bmax q\equ p\bmax|p-q|$.
\end{proof}

We state our results first for the simpler case of sequences of
independent indicator variables with given (non-random) probabilities.
The following theorem gives necessary and sufficient conditions for
asymptotical equivalence and contiguity. 
(The asymptotical equivalence criterion follows by a simple and
standard type of 
calculation with Hellinger distances, see the proof in \refS{Spf}
and, \eg, \cite[p.\ 158]{OZ}, although we have not seen it stated in
this form before.
The contiguity criterion is a special case of a result by \citet{OZ}
for general sequences of independent variables.)
The proofs of the theorems are given in \refS{Spf}.

\begin{theorem}
  \label{T1}
Let $1\le N(n)\le\infty$ and let $X_n=(\ini)\xiinn$ and 
$X'_n=(\ini')\xiinn$ be finite or infinite random vectors 
consisting of independent indicator variables $\ini\sim\Be(\pni)$ and
$\ini'\sim\Be(\pni')$.
\begin{romenumerate}
\item
  $X_n\aseq X_n'$ if and only if
\begin{equation}
	\label{t1i}
\suminn\rho(\pni,\pni')\to0.
  \end{equation}
\item
$X_n\cont X_n'$ if and only if
\begin{equation}
	\label{t1iia}
\suminn\rho(\pni,\pni')=O(1)
  \end{equation}
and, with $\qni\=1-\pni$ and $\qni'\=1-\pni'$,
\begin{equation}
	\label{t1iib}
\lim_{C\to\infty}
\limsup_\ntoo
\lrpar{
\sum_{i:\pni>C\pni'}\pni
+\sum_{i:\qni>C\qni'}\qni}
=0.
  \end{equation}
\end{romenumerate}
\end{theorem}

\begin{remark}
  By symmetry, $X_n\contr X_n'$ is equivalent to \eqref{t1iia} and
\begin{equation}
	\label{t1iic}
\lim_{C\to\infty}
\limsup_\ntoo
\lrpar{
\sum_{i:\pni'>C\pni}\pni'
+\sum_{i:\qni'>C\qni}\qni'}
=0,
  \end{equation}
and thus $X_n\contxx X_n'$ is equivalent to \eqref{t1iia},
\eqref{t1iib} and \eqref{t1iic}.
\end{remark}

\begin{remark}
  Often $\pni\le0.9$ for all $n$ and $i$, and then the second sum in
  \eqref{t1iib} vanishes for $C>10$ and can thus be omitted.
\end{remark}

\begin{remark}
  The condition \eqref{t1iib} is only needed to take care of cases
  when $\pni$ and $\pni'$  (or $\qni$ and $\qni'$, in case $\pni$ and
  $\pni'$ are close to 1) are not of the same order. If no such $\pni$
  and $\pni'$ appear, which is the typical case, then \eqref{t1iia} is
  thus enough. 

We may rewrite \eqref{t1iib} in several ways. 
For example, it is equivalent to (following the formulation in
\cite{OZ} in a more general case):
for every sequence $C_n\to\infty$,
\begin{equation}
	\label{t1iiboz}
\sum_{i:\pni>C_n\pni'}\pni
+\sum_{i:\qni>C_n\qni'}\qni
\to0.
  \end{equation}
It is also equivalent to:
For every $\eps>0$, there exist $C$ and $n_0$ such that if $n\ge n_0$, then
\begin{equation}
	\label{t1iibr}
\sum_{i:\pni>C\pni'}\pni
<\eps
\qquad
\text{and}
\qquad
\sum_{i:\qni>C\qni'}\qni
< \eps.
  \end{equation}
\end{remark}

\begin{remark}
  As pointed out by \citet{OZ}, \eqref{t1iia} does not imply
  \eqref{t1iib} in general. A simple counter example is provided by 
$N(n)=n$, $\pni=n\qw$, $\pni'=n\qww$.
(On the other hand, it is easy to see, and also follows by the theorem,
  that \eqref{t1i} implies \eqref{t1iib} and \eqref{t1iic}.)
\end{remark}

\begin{remark}\label{Rkakutani}
  In the very special case when $N(n)$, $\pni$ and $\pni'$ do not
  depend on $n$ (and we omit the subscript $n$), it is easily shown
  that if $0<p_i'<1$ for all $i$, then \eqref{t1iia} implies
  \eqref{t1iib}, and thus $(X)\cont(X')$, which by \refE{EA1} says
  that the distribution of $X$ is absolutely continuous with respect
  to the distribution of $X'$.
If we further assume also $0<p_i<1$, by symmetry the distributions
  are thus mutually absolutely continuous.
This is part of the dichotomy by Kakutani for product measures,
see \eg{} \cite[Corollary IV.2.38]{JS},
which in our case says that either 
$\sum_i\rho(p_i,p_i')<\infty$ and
  the distributions are mutually absolutely continuous, or
$\sum_i\rho(p_i,p_i')=\infty$ and
  the distributions are mutually singular.

Returning to the general case in \refT{T1}, it is easy to show that,
analoguously, if $\sum_i\rho(\pni,\pni')\to\infty$, then the
distributions of $X_n$ and $X_n'$ are asymptotically mutually singular
in the sense that there exist sets $A_n$ with $\P(X_n\in A_n)\to1$ and
$\P(X_n'\in A_n)\to0$, \cf{} \cite[Theorem V.2.32]{JS}.
\end{remark}

\begin{remark}
  We have stated  \refT{T1} in terms of sequences of
pairs of random vectors. It is possible (at least partly) to rephrase
it in terms of estimates for a single pair $(X,X')$, see Lemmas
\refand{L1i}{L1ii} below. Similar reformulations may be made for
\refT{T2}, but we leave these to the reader.
\end{remark}

We extend \refT{T1} to the case of random probabilities $\pni$. In
this case we cannot expect conditions that are both necessary and
sufficient, so we give only sufficient conditions, which are more
important in applications.
(An important obstacle to finding necessary conditions is that
different distributions of the probabilities may give the same
distribution of the 
indicators. For example, using the notation of \refT{T2}, if $\pni$
are \iid{} with $\pni\sim U(0,1)$ and $\pni'=1/2$, then $X_n\eqd X_n'$.)

\begin{theorem}
  \label{T2}
Let $1\le N(n)\le\infty$ and 
suppose that $\pp_n=\pnix$ and $\pp_n'=\pnixy$ are random vectors in
$\oi^{N(n)}$. 
Let $X_n=(\ini)\xiinn$ and 
$X'_n=(\ini')\xiinn$ be random vectors 
of indicator variables such that the conditioned random vectors
$(X_n\mid\pp_n)$ and $(X_n'\mid\pp_n')$ are sequences of 
independent indicator variables with 
$(\ini\mid\pp_n)\sim\Be(\pni)$ and
$(\ini'\mid\pp_n')\sim\Be(\pni')$.
\begin{romenumerate}
\item
If
\begin{equation}
	\label{t2i}
\suminn\rho(\pni,\pni')=o_p(1),
  \end{equation}
then  $X_n\aseq X_n'$.
\item
If
\begin{equation}
	\label{t2iia}
\suminn\rho(\pni,\pni')=O_p(1)
  \end{equation}
and, with $\qni\=1-\pni$ and $\qni'\=1-\pni'$, for every $\eps>0$,
\begin{equation}
	\label{t2iib}
\lim_{C\to\infty}
\limsup_\ntoo
\P\lrpar{
\sum_{i:\pni>C\pni'}\pni
+\sum_{i:\qni>C\qni'}\qni>\eps}
=0,
  \end{equation}
then $X_n\cont X_n'$.
\end{romenumerate}
\end{theorem}

\begin{remark}\label{Rop}
Recall 
that if $S_n$ denotes the random sum on the \lhs{} of 
\eqref{t2i}, then  \eqref{t2i} can also be written
$S_n\pto0$.
Similarly, 
the $O_p(1)$ notation in \eqref{t2iia} 
means that for every $\eps>0$, there exists
  $C$ such that $\P(S_n>C)<\eps$ for all $n$; this is also known as
  stochastic boundedness or tightness of the sequence \set{S_n}, and
  is equivalent to $\P(S_n>C_n)\to0$ for every sequence $C_n\to\infty$.
\end{remark}

\begin{remark}
In analogy to \eqref{t1iiboz}, the condition \eqref{t2iib} is
equivalent to:
For every sequence $C_n\to\infty$,
  \begin{equation}
	\label{t2iiboz}
\sum_{i:\pni>C_n\pni'}\pni
+\sum_{i:\qni>C_n\qni'}\qni
\pto0.
  \end{equation}
\end{remark}

As said above, Theorems \refand{T1}{T2} apply immediately to random
graphs $\gnpp$. We state a version of \refT{T2} for this case, where
we have added some simplifying assumptions.
Recall that $\pij$ and $\pij'$ may (and typically do) depend on $n$,
although we do not show that in our notation. 

\begin{corollary}
  \label{C1}
Let, for each $n$, $\pp=\pijnx$ and $\pp'=\pijnxy$ be random vectors
of probabilities
and suppose that \whp{} $\max_{i,j}\pijn\le0.9$.
\begin{romenumerate}
\item
  If 
  \begin{equation}\label{c1i}
\sum_{i<j} \frac{(\pijn-\pijn')^2}{\pijn} =o_p(1),	
  \end{equation}
then $\gnpp\aseq\gnppy$.
\item
  If 
  \begin{equation}\label{c1ii}
\sum_{i<j} \frac{(\pijn-\pijn')^2}{\pijn} =O_p(1),	
  \end{equation}
then $\gnpp\contr\gnppy$.
\item
If \eqref{c1ii} holds, and further, for some constant $c>0$,
\whp{} $c\pij\le \pij' \le 0.9$ for all $i,j$, 
then $\gnpp\contxx\gnppy$.
\end{romenumerate}
\end{corollary}

We specialize further to an important case.

\begin{corollary}
  \label{C2}
Let, for each $n$, $\pp=\pijnx$ and $\pp'=\pijnxy$ be random vectors
of probabilities
and suppose that $\pijny=\pijn+O(\pijn^2)$.
\begin{romenumerate}
\item
  If\/ $\sum_{i<j} \pijn^3 =o_p(1)$,	
then $\gnpp\aseq\gnppy$.
\item
  If\/ $\sum_{i<j} \pijn^3 =O_p(1)$,	
and further, for some constant $c>0$,
\whp{} $\max_{i,j}\pij\le0.9$, $\max_{i,j}\pij'\le0.9$ and
$\pij'\ge c\pij$ for all $i,j$,
then $\gnpp\contxx\gnppy$.
\end{romenumerate}
\end{corollary}

\section{Examples}\label{Sex}

\begin{example}
  \label{Ekappa}
\citet{SJ178} study a general class of sparse random graphs $G(n,\gk)$
which include many cases studied earlier by various authors. These random
graphs are 
defined as \gnpij{} with
\begin{equation}
  \label{pij1}
\pij\=\pijq1\=\min\Bigpar{\frac{\gk(x_i,x_j)}{n},1}
=\pijo\bmin1,
\end{equation}
with \begin{equation}
  \label{pij0}
\pijo\=\frac{\gk(x_i,x_j)}{n},
\end{equation}
where $\gk:\cSS\times\cSS\to\ooo$ is a symmetric measurable function
defined on some measurable space $\cSS$ and $x_1,\dots,x_n$ is a
random sequence of elements of $\cSS$, not necessarily \iid{} but such
that the empirical distribution of $x_1,\dots,x_n$ converges to a
probability measure $\mu$ on $\cSS$; see \cite{SJ178} for details. 
(Some further technical conditions are assumed in \cite{SJ178}; they
are not relevant here.)
Typically, \whp{} $\gk(x_i,x_j)\le n$ for all $i,j$, and then
$\pij$ equals the simpler $\pijo$.
Two natural variations, also treated in \cite{SJ178} and used in 
various cases by various authors,
are obtained by replacing \eqref{pij1} by
\begin{align}\label{pij2}
  \pijq2&\=1-\exp\Bigpar{-\frac{\gk(x_i,x_j)}{n}}
=1-\exp\bigpar{-\pijo}
\intertext{or}\label{pij3}
\pijq3&\=\frac{\gk(x_i,x_j)}{n+\gk(x_,x_j)}
=\frac{\pijo}{1+\pijo}.
\end{align}
(Thus $\pijq3/(1-\pijq3)=\pijo$; at least in the case studied in
\refE{Erank1} below, this is in some sense
simpler and more natural, see \citet{BrittonDML}.)
In all cases $\pijq\ell=\pijo+O(\pijo^2)$, which is the essential
estimate for our purposes;
the results below extend to the general case
\begin{equation}\label{pijh}
  \pij\=h(\pijo)
\quad\text{for a function $h$ with}\quad
h(x)=x+O(x^2).
\end{equation}

It was shown in \cite{SJ178} that the same asymptotic results hold
for these three versions for the properties studied there. We can now
show that, under an extra condition, the three versions 
are \aseqq, and thus have the
same asymptotic behaviour for \emph{any} property.
Indeed, \refC{C2} applies immediately and shows that if 
\begin{equation}
  \label{gk3o}
\sum_{1\le i<j\le n} \gk(x_i,x_j)^3 = o_p(n^3),
\end{equation}
then all three $G(n,\pijq\ell)$, $\ell=1,2,3$, are \aseqq;
similarly, if the weaker 
\begin{equation}
  \label{gk3O}
\sum_{1\le i<j\le n} \gk(x_i,x_j)^3 = O_p(n^3)
\end{equation}
holds together with $\max_{i,j}\pijo\le0.9$ \whp,
then all three $G(n,\pijq\ell)$ are mutually contiguous.

In fact, \eqref{gk3O} alone suffices for 
$G(n,\pijq2)\contxx G(n,\pijq3)$ 
because \eqref{gk3O} implies 
$\max_{i,j}\pijo=O_p(1)$ so by conditioning we may assume that
$\max_{i,j}\pijo\le C_1$ for some constant $C_1$, and then, for
$\ell=2,3$, $\pijq\ell\le C_2<1$ and $c\pijo\le\pijq\ell\le\pijo$.
Furthermore, by the same conditioning and \refC{C1}(ii), \eqref{gk3O}
implies 
$G(n,\pijq2)\contr G(n,\pijq1)$  and $G(n,\pijq3)\contr G(n,\pijq1)$.
However, if, for example $\gk(x_1,x_2) \ge n$ \whp, then 
$I_{12}=1$ \whp{} in  $G(n,\pijq1)$ but not in $G(n,\pijq2)$ or $G(n,\pijq3)$,
and we do not have contiguity in the opposite direction.

We study some special cases in the following examples.
\end{example}

\begin{example}\label{Eiid}
One common case of the construction in \refE{Ekappa}
uses $x_1,\dots,x_n$ that are
\iid{} on $\cSS$ with distribution $\mu$. In this
case, we show that the condition
\begin{equation}
  \label{gkt2o}
\mu\times\mu\set{(x,y):\gk(x,y)>t)}=o\bigpar{t\qww} \qquad\asttoo
\end{equation}
implies \eqref{gk3o} and thus asymptotic equivalence of the three
versions.
In particular, this holds if $\int_{\cSS\times\cSS}
\gk(x,y)^2\,d\mu(x)\,d\mu(y)<\infty$.

In fact, if $G(t)\=\mu\times\mu\set{(x,y):\gk(x,y)>t)}=o(t\qww)$, then
\begin{equation}\label{magnus}
  \begin{split}
\E \bigpar{\gk(x_1,x_2)\bmin n}^3
&=\int_{\cSS\times\cSS} \bigpar{\gk(x,y)\bmin n}^3\,d\mu(x)\,d\mu(y)
\\&
=\int_0^n 3 t^2 G(t) \,dt 	
=n\int_0^1 3 (ns)^2 G(ns) \,ds
=o(n) 	
  \end{split}
\end{equation}
by \eqref{gkt2o} and dominated convergence. Hence,
$\E\sum_{i<j} (\gk(x_i,x_j)\bmin n)^3 =o(n^3)$, so 
$\sum_{i<j} (\gk(x_i,x_j)\bmin n)^3 =o_p(n^3)$. Moreover,
\begin{multline*}
\P\Bigpar{\sum_\ijn (\gk(x_i,x_j)\bmin n)^3 
\neq \sum_{\ijn} \gk(x_i,x_j)^3} 
\\
\le\sum_{\ijn}\P(\gk(x_i,x_j)> n)
\le n^2G(n)=o(1),
\end{multline*}
and \eqref{gk3o} follows.

Similarly, we can easily shown that \eqref{gk3O}, and thus at least
partial contiguity, follows from 
\begin{equation}
  \label{gkt2O}
\mu\times\mu\set{(x,y):\gk(x,y)>t)}=O\bigpar{t\qww} \qquad\asttoo.
\end{equation}
In this case, given $\eps>0$, there exists $C_1$ such that 
$$\P\Bigpar{\sum_\ijn (\gk(x_i,x_j)\bmin C_1n)^3 \neq \sum_{\ijn}
  \gk(x_i,x_j)^3}  \le\eps;$$ 
further, a calculation as in
\eqref{magnus}
yields
$\E\sum_\ijn (\gk(x_i,x_j)\bmin C_1n)^3 =O(n^3)$
and thus
$\P\bigpar{\sum_\ijn (\gk(x_i,x_j)\bmin C_1n)^3>C_2 n^3}<\eps$ for
  some $C_2$; we omit the details.
\end{example}

\begin{example}\label{Edet}
  Another case of the construction in \refE{Ekappa}
uses $\cSS=(0,1]$ with $\mu$ = Lebesgue measure and the deterministic
$x_i=i/n$, $i=1,\dots,n$. 
The homogeneous case $\gk(x,y)=c/(x\bmax y)$ yielding
$\pijo=c/(i\bmax j)$, where $c>0$ is a constant, is particularly
interesting
and related to the CHKNS model,
see \citet[Sections 16.1]{SJ178},
\citet{Durrett03,Durrett}
and 
\citet{Rsmall} 
and the
references given there.

In this case, $\sum_{1\le i<j<\infty} \pij^3 \le
c^3\sum_{j\ge2} j\cdot j^{-3}<\infty$, and thus
$\sum_{ i<j} \pij^3 =O(1)$;
if we further for simplicity assume
$c<2$ and thus $\max_{ij}\pijo\allowbreak<1$, then \refC{C2}(iii)
implies that 
$G(n,\pijq1)\contxx G(n,\pijq2)\contxx G(n,\pijq3)$.

Note that in this case, $p_{12}\supx1$, $p_{12}\supx2$ and $p_{12}\supx3$
are constant and different, which shows that the three random graphs
are \emph{not} \aseqq{} (for a trivial reason).

We have $\pijq3=c/(i \bmax j +c)$; the same results hold for the
further variation $\pij=c/(i \bmax j + d)$ for any $d>c-2$.

In this example, the infinite random graphs $G(\infty,\pijq\ell)$,
$\ell=1,2,3$, 
are well-defined too, and it follows from Kakutani's criterion
discussed in \refR{Rkakutani} that (still provided $c<2$)
these three infinite random graphs have
mutually absolutely continuous distributions,
which is the infinite graph version of the contiguity result just
given for finite $n$, \cf{} \refE{EA1}.
(The infinite random graph $G(\infty,\pijq1)$ was studied before the
finite version, see \cite{KW} and \cite{Durrett03},
\cite{Durrett} with further references.)
\end{example}

\begin{example}\label{Edet2}
A related case 
uses the same $\cSS=(0,1]$, $\mu$ = Lebesgue measure and 
$x_i=i/n$, $i=1,\dots,n$, as \refE{Edet}, now with
the homogeneous 
$\gk(x,y)=c/\sqrt{xy}$ yielding
$\pijo=c/\sqrt{ij}$; this case is a mean-field version of
the preferential
attachment model by Barab\'asi and Albert
\cite{BAsc},
see \citet[16.2]{SJ178} and \citet{Rsmall} 
and the
references given there.

Also in this case, $\sum_{1\le i<j<\infty} \pij^3 <\infty$, and thus
$\sum_{ i<j} \pij^3 =O(1)$ (in spite of the fact that \eqref{gkt2O}
does not hold);
if we further for simplicity assume
$c<\sqrt2$, and thus $\max_{ij}\pijo\allowbreak<1$, we obtain 
the same results as in \refE{Edet}.
\end{example}

\begin{example}\label{Erank1}
A common case of \refE{Ekappa}
is when $\gk(x,y)=\psi(x)\psi(y)$ for some function
$\psi:\cSS\to\ooo$, see \cite[Section 16.4]{SJ178} for discussion and
references to previous papers.

In this case, $\sum_{i<j}\gk(x_i,x_j)^3 \le
\bigpar{\sum_i\psi(x_i)^3}^2$,
so \eqref{gk3o} and \eqref{gk3O} may be replaced by
\begin{align}\label{psi3o}
  \sumin \psi(x_i)^3&=o_p\bigpar{n^{3/2}},
\intertext{and}
\label{psi3O}
  \sumin \psi(x_i)^3&=O_p\bigpar{n^{3/2}}.
\end{align}

If we combine this choice of $\gk$ with the \iid{} choice of $x_i$ in
\refE{Eiid}, 
it is easily seen, arguing as in \eqref{magnus} but now with
$\sum_i\bigpar{\psi(x_i)\bmin n\qq}$, that
\begin{equation}
  \label{gpsi2o}
\mu\set{x:\psi(x)>t)}=o\bigpar{t\qww} \qquad\asttoo
\end{equation}
implies \eqref{gk3o} and thus asymptotic equivalence of the three
versions; in particular this holds if 
$\int \psi(x)^2\,d\mu(x)<\infty$.
Similarly,
\begin{equation}
  \label{gpsi2O}
\mu\set{x:\psi(x)>t)}=O\bigpar{t\qww} \qquad\asttoo
\end{equation}
implies \eqref{gk3O} and thus at least partial contiguity.
\end{example}

\begin{example}
\label{EEHH}
\citet{EHH} study a minor variation of the construction in
\refE{Erank1}; they let $\gL_1,\dots,\gL_n$ be positive \iid{} random
variables
with some fixed distribution 
and define $\pij$ by (in our notation) 
\eqref{pij1}, \eqref{pij2},  \eqref{pij3} or more generally
\eqref{pijh} with 
\begin{equation}\label{pijgL}
  \pijo\=\frac{\gL_i\gL_j}{\sum_1^n \gL_i}.
\end{equation}
(This too can be seen as an instance of the general construction in
\refE{Ekappa}, see \cite[Section 16.4]{SJ178}.)

Assume that $\P(\gL_1>t)=o(t\qww)$ (which is the case in \cite{EHH}). 
Then, just as \eqref{psi3o}
follows from \eqref{gpsi2o}, $\sum_1^n \gL_i^3=o_p(n^{3/2})$. Since
further 
$\sum_1^n \gL_i/n\pto \E\gL>0$ by the law of large numbers,
it follows from \eqref{pijgL} that $\sum_{i<j}\pij^3\pto0$. 
Hence, if we compare two random graphs $\gnpp$ and $\gnppy$ defined by
this method for two different functions $h$ and $h'$, we obtain 
$\gnpp\aseq\gnppy$
by \refC{C2}.

\citet{EHH} study the distance $H_n$ between two random points, and (a
minor) part of their proof consists in showing that the choice of $h$ does not
matter (enabling them to consider only the version \eqref{pij2} in the
main part of the proof):
the variables $H_n$ and $H_n'$ obtained by two different functions $h$
and $h'$ in \eqref{pijh} can be coupled such that $\P(H_n\neq
H_n')=o(1)$, or in our notation $H_n\aseq H_n'$, see \refT{T0}.
We thus obtain this as an immediate consequence of the stronger
statement $\gnpp\aseq\gnppy$, which by \refT{T0}
means that the random graphs can
be coupled with $\P\bigpar{\gnpp\neq\gnppy}\to0$.
\end{example}

\begin{example}
  Our results are stated for graphs with a deterministic number of
  vertices, but can be extended to graphs with random vertex set too
  by conditioning on the vertex set. One interesting
such case is obtained from
  \refE{Ekappa} by letting $x_1,\dots,x_n$ be the points of a Poisson
  process on $\cS$ with intensity $\gl\mu$, where $\gl>0$ is our
  parameter and we consider asymptotics as $\gl\to\infty$; thus $n$ is
  random with the distribution $\Po(\gl)$.

Conditioned on $n$, we have the situation in \refE{Eiid}. It follows,
for example, that if \eqref{gkt2o} holds, then the random graphs
defined in this way using \eqref{pij1}, \eqref{pij2} and \eqref{pij3}
are \aseqq; we omit the details.
\end{example}

\begin{example}
  \label{Eclustering}
\citet{clustering} study a generalization of the model in
\refE{Ekappa} where small sets of edges are added at once, thus
allowing a certain degree of clustering; more precisely, for every
subgraph $F$ of the complete graph $K_n$, we have a certain
probability of adding (the edges of) $F$, and these events are
independent for different $F$.
While this introduces dependencies between the edge indicators,
the results of the present paper are still applicable to the sequence
of indicators $I_F$ describing the added sets of edges, and 
asymptotic equivalence
or contiguity for two versions of this sequence obviously implies 
asymptotic equivalence or contiguity for the resulting random graphs too.

We leave the explicit statement of results in this case to the reader.
\end{example}

\begin{example}
\label{Ep2}
In this final example, let us return to the case of deterministic
$\pp=\pijnx$ and let us change all $\pij$ proportionately to
$\pij'\=(1+\gdn)\pij$ for some $\gdn$. Assume for simplicity that all
$\pij\le0.9$ and that $|\gdn|\le1$.

By \refC{C1}(i), if 
$\gdn^2\sum_{i<j}\pij\to0$, then
$\gnpp\aseq\gnppy$.
Further, by \refC{C1}(iii), if 
$\gdn^2\sum_{i<j}\pij=O(1)$ and, for simplicity, $\gdn\to0$, then
$\gnpp\contxx\gnppy$.

In fact, by \eqref{rho11}, 
$\rho(\pij,\pij')\equ \gdn^2\pij \bmin |\gd_n|\pij = \gdn^2\pij$, and
thus 
 by \refT{T1} 
the conditions
$\gdn^2\sum_{i<j}\pij\to0$
and 
$\gdn^2\sum_{i<j}\pij=O(1)$
are necessary too for asymptotic equivalence and contiguity,
respectively.  
(The necessity can also be checked by considering the total number of
edges, as in the special case in \refE{Egnp}.)
\end{example}

\begin{remark}
As in \refE{Ep2}, necessity in \refT{T1} can in many cases where
$\pij'\le\pij$ for all $i$ and $j$ (or conversely)
be proved by considering
the total numbers $\sum_i\ini$ and $\sum_i\ini'$, but this method does
not suffice in all cases. A simple counter example is given by $N(n)=n^4+n^8$,
$\pni=n\qw$ for $1\le i\le n^4$ and $\pni=n^{-3}$ for $i>n^4$, and
$\pni'=\pni-\pni^2$; it is easily checked that 
then \eqref{t1i} and \eqref{t1iia} do not hold, and thus we do not
have asymptotic equivalence or even contiguity, but, using 
\cite[Theorems 2.M and 1.C]{SJI},
\begin{equation*}
  \dtv\Bigpar{\sum_i\ini,\sum_i\ini'}
=\dtv\bigpar{\Po(n^5+n^3),\Po(n^5+n^3-2n^2)}+o(1)
\to0.
\end{equation*}

\end{remark}

\section{More on asymptotic equivalence and contiguity}
\label{Smore}

We will use two metrics to measure the distance between probability
distributions. We state some well-known definitions and facts, see
\eg{} \cite[Appendix A.1]{SJI} and \cite[Chapter IV.1 and V.4a]{JS}.

\begin{definition}\label{D2}
  If $P$ and $Q$ are two probability measures on the same measurable
  space $(\cXX,\cA)$, and $R$ is any $\gs$-finite measure on
  $(\cXX,\cA)$ such that $P\ll R$ and $Q\ll R$, define
the \emph{total variation distance}
\begin{equation}\label{dtv}
  \dtv(P,Q)\=\sup_{A\in\cA}|P(A)-Q(A)|
=\frac12\int_{\cXX}\lrabs{\frac{dP}{dR}-\frac{dQ}{dR}}dR
\end{equation}
and the \emph{Hellinger distance}
\begin{equation}\label{dh}
  \dh(P,Q)\=\lrpar{\frac12\int_{\cXX}
 \lrpar{\sqrt{\frac{dP}{dR}}-\sqrt{\frac{dQ}{dR}}}^2dR}\qq
=\bigpar{1-H(P,Q)}\qq
\end{equation}
where $H(P,Q)$ is the \emph{Hellinger integral}
\begin{equation}\label{H}
  H(P,Q)\=\int_{\cXX}\sqrt{\frac{dP}{dR}}\sqrt{\frac{dQ}{dR}}\,dR.
\end{equation}
(We can, at least symbolically,  write \eqref{H} as
$H(P,Q)\=\int_{\cXX}\sqrt{dP\,dQ}$.)
Note that these quantities do not depend on the choice of $R$. (We may
thus take, \eg, $R=P+Q$.)

We have $\dtv(P,Q)=\frac12\norm{P-Q}$, using the standard norm on
real-valued measures. (The factor $\frac12$ is conventional and
convenient but unimportant, as is the factor $\frac12$ in the definition
of $\dh$.)

We use the same notations for two random variables $X$ and $Y$ with
values in $\cXX$, with $\dtv(X,Y)\=\dtv(\cL(X),\cL(Y))$ and similarly
for $\dh$ and $H$. (Thus $\dtv(X,Y)=0\iff\dh(X,Y)=0\iff X\eqd Y$.)
In particular,
\begin{equation}\label{dtvxy}
  \dtv(X,Y)\=\sup_{A\in\cA}|\P(X\in A)-\P(Y\in A)|.
\end{equation}
\end{definition}

It is easily seen that
$\dtv$ and $\dh$ are metrics on the set of all probability measures
on $(\cXX,\cA)$; further, $0\le\dtv\le1$, $0\le\dh\le1$ and $0\le H \le1$,
and 
\begin{equation}\label{dtvh}
  \dhh(P,Q)\le\dtv(P,Q)\le\sqrt2\,\dh(P,Q);
\end{equation}
hence $\dtv$ and $\dh$ are equivalent metrics.
Furthermore, $\dtv(P,Q)=1\iff\dh(P,Q)=1\iff P\perp Q$, \ie, $P$ and $Q$ are
mutually singular.

Recall that a \emph{coupling} of two random variables $X$ and $Y$ with
values in the same space is a pair of random variables $(X',Y')$,
defined together on the same probability space, with $X'\eqd X$ and $Y'\eqd Y$.
It is well-known that
\begin{equation}\label{coupling}
  \dtv(X,Y)=\min_{(X',Y')}\P(X'\neq Y'),
\end{equation}
taking the minimum over all couplings $(X',Y')$ of $X$ and $Y$.

\begin{theorem}\label{T0}
  Let\/ $X_n$ and\/ $Y_n$ be random variables with values in $\cXX_n$. Then
  the following are equivalent.
  \begin{romenumerate}
\item\label{T0a}
$X_n\aseq Y_n$.
\item\label{T0b}
$\dtv(X_n,Y_n)\to0$.
\item\label{T0c}
$\dh(X_n,Y_n)\to0$.
\item\label{T0e}
$H(X_n,Y_n)\to1$.
\item\label{T0d}
There exist couplings $(X_n',Y_n')$ of $X_n$ and
$Y_n$ such that $\P(X'_n\neq Y'_n)\to0$. 
  \end{romenumerate}
\end{theorem}
\begin{proof}
This too is well-known and easy:
  \ref{T0a}$\iff$\ref{T0b} by \eqref{dtvxy} and \refD{D1};
  \ref{T0b}$\iff$\ref{T0c} by \eqref{dtvh};
  \ref{T0c}$\iff$\ref{T0e} by \eqref{dh};
  \ref{T0b}$\iff$\ref{T0d} by \eqref{coupling}.
\end{proof}

We calculate
the Hellinger distance and integral for two Bernoulli
distributions. (This is the origin of our function $\rho$ in \refD{Drho}.)
\begin{lemma}\label{LBe}
For any $p,q\in\oi$,
  \begin{align*}
  \dh\bigpar{\Be(p),\Be(q)}
&=2\qqw\rho(p,q)\qq,
\\
  H\bigpar{\Be(p),\Be(q)}
&=1-\tfrac12\rho(p,q).
  \end{align*}
\end{lemma}
\begin{proof}
  Use \eqref{dh} with $P=\Be(p)=p\gd_0+(1-p)\gd_1$, 
$Q=\Be(q)=q\gd_0+(1-q)\gd_1$ and $R=\gd_0+\gd_1$, together with the
  definition \eqref{rho1}.
\end{proof}

An important, and well-known, property of Hellinger distances and
integrals is that they behave simple for product measures.
Let
$[n]\=\set{1,\dots,n}$ when $n<\infty$ and $[\infty]\=\bbN=\set{1,2,\dots}$.

\begin{lemma}\label{Lprod}
Let $1\le N\le\infty$ and let, for $i\in[N]$,
$P_i$ and $Q_i$ be  probability measures on the same
measurable space $(\cXX_i,\cA_i)$.
If $P=\prodiN P_i$ and $Q=\prodiN Q_i$, then
$
H(P,Q)=\prodiN H(P_i,Q_i)$.
\end{lemma}

\begin{proof}
  This is stated in, \eg, \cite[Proposition IV.1.73]{JS}, but for
  completeness we give the simple proof.

If $N<\infty$, the result is an immediate consequence of \eqref{H} and
Fubini's theorem, choosing \eg{} $R_i=P_i+Q_i$ and $R=\prodiN R_i$.

If $N=\infty$, let $\cF_n$ be the $\gs$-field on $\prodoo \cXX_i$ given
by \set{A\times\prod_{n+1}^\infty\cXX_i:A\in\prod_1^n\cA_i}, and let
$\bp_n\=P|_{\cF_n}$ and 
$\bq_n\=Q|_{\cF_n}$.
Then, using the finite case,
\begin{equation*}
 H(\bp_n,\bq_n)=H\Bigpar{\prodn P_i,\prodn Q_i}
=\prodn H(P_i,Q_i). 
\end{equation*}
Furthermore, choosing $R=(P+Q)/2$ on $\cXX\=\prodoo\cXX_i$,
$d\bp_n/dR=\E(d P/dR\mid\cF_n)$ with respect to $R$, so $(d\bp_n/dR)$
is a bounded $R$-martingale and $d\bp_n/dR\to dP/dR$ $R$-\as, and
similarly for $\bq_n$. Hence, \eqref{H} and dominated convergence
yields $H(P,Q)=\lim_\ntoo H(\bp_n,\bq_n)=\prodoo H(P_i,Q_i)$. 
\end{proof}

\section{Proofs}\label{Spf}

\begin{proof}[Proof of \refT{T1}]
\pfitem{i}
By Lemmas \refand{Lprod}{LBe},
\begin{equation*}
  H(X_n,X_n')=\prodNn H(\ini,\ini')
=\prodNn \bigpar{1-\tfrac12\rho(\pni,\pni')}.
\end{equation*}
Hence, 
\begin{equation*}
1-\tfrac12 \sumNn \rho(\pni,\pni')
\le
  H(X_n,X_n')
\le
\exp\biggpar{-\tfrac12 \sumNn \rho(\pni,\pni')},
\end{equation*}
and thus $H(X_n,X_n')\to1 \iff \sumNn \rho(\pni,\pni')\to0$, which
yields the result by \refT{T0}.

\pfitem{ii}
This is, in view of \refL{LBe} and the equivalence of \eqref{t1iib}
and \eqref{t1iiboz},
a special case of \cite[Theorem 1]{OZ}, to which we refer for
a complete proof.
Nevertheless, for completeness, we sketch a proof of the more important
``if'' direction.

First, we can by a simpler version of the argument in the proof of
\refT{T2} below assume that $\pni\le C_2\pni'$ and $\qni\le C_2\qni'$
for some constant $C_2$. (We define $\pni'''$ by \eqref{p'''} with
$\pni''\=\pni$ and use \eqref{ew}--\eqref{emma}.)
Under this assumption, if we let 
$\Pni\=\cL(\ini)=\Be(\pni)$, 
$\Pn\=\prod_i\Pni$,
$\Pni'\=\cL(\ini')=\Be(\pni')$, 
$\Pn'\=\prod_i\Pni'$, 
we have by Fubini, using $\int (d\Pni/d\Pni')\,d\Pni' =1$ and \eqref{rho2},
\begin{align*}
  \int \lrpar{\frac{d\Pn}{d\Pn'}}^2 d\Pn'
&=
\prod_i  \int \biggpar{\frac{d\Pni}{d\Pni'}}^2 d\Pni'
=
\prod_i \lrpar{1+ \int \biggpar{\frac{d\Pni}{d\Pni'}-1}^2 d\Pni'}
\\&
=
\prod_i 
\biggpar{1+ \parfrac{\pni-\pni'}{\pni'}^2\pni' 
 +\parfrac{\qni-\qni'}{\qni'}^2\qni'}
\\&
=
\prod_i 
\biggpar{1+\frac{(\pni-\pni')^2}{\pni'} +\frac{(\pni-\pni')^2}{1-\pni'} }
\\&
\le
\prod_i 
\biggpar{1+(C_2+1)\frac{(\pni-\pni')^2}{\pni+\pni'} 
+(C_2+1)\frac{(\pni-\pni')^2}{1-\pni+1-\pni'} }
\\&
\le
\prod_i 
\bigpar{1+(C_2+1)C\rho(\pni,\pni')}
\\&
\le
\exp\Bigpar{(C_2+1)C\sum_i \rho(\pni,\pni')}
\\&
=O(1)
\end{align*}
and thus for any sets $A_n$, by the \CSineq,
\begin{align*}
P_n(A_n)
&=\int_{A_n} d\Pn
\le\lrpar{\int\parfrac{d\Pn}{d\Pn'}^2 d\Pn'\cdot\int_{A_n}d\Pn'}\qq
=O\bigpar{\Pn'(A_n)\qq}
\end{align*}
and thus $P_n\cont P_n'$, which is the same as
$X_n\cont X_n'$.
\end{proof}

We say that a finite or infinite
random vectors of indicator variables 
$X=(I_i)\xiin$ has distribution
$\Bex(\pp)$, where 
$\pp=\set{p_i}\xiin$ is a
deterministic vector with elements in $\oi$, 
if the random variables
$I_i$ are
independent indicator variables with 
$I_i\sim\Be(p_i)$.

More generally, 
if $\pp=\set{p_i}\xiin$ is a
random vector with elements in $\oi$, 
with $N\le\infty$,
we say that 
random vectors of indicator variables 
$X=(I_i)\xiin$ has distribution
$\Bex(\pp)$
if the conditioned random vector
$(X\mid\pp)$ is a sequence of 
independent indicator variables with 
$(I_i\mid\pp)\sim\Be(p_i)$.

We next give two results comparing two random vectors with
distributions $\Bex(\pp)$ and $\Bex(\pp')$ with deterministic $\pp$
and $\pp'$. The first result is easily seen to be equivalent to
the ``if'' direction of
\refT{T1}(i), while the second is equivalent to a special case of
the ``if'' direction of
\refT{T1}(ii).

\begin{lemma}
  \label{L1i}
For every $\eps>0$, there exists $\gd>0$ such
that if $X\sim\Bex(\pp)$ and $X'\sim\Bex(\pp')$ for two deterministic
vectors $\pp=\set{p_i}\xiin$ and $\pp'=\set{p_i'}\xiin$ of the same
length $N\le\infty$, and these satisfy
$
\sumiN\rho(p_i,p_i')< \gd
$,
then 
$\dtv(X,X')<\eps$.
\end{lemma}
\begin{proof}
  Suppose not. Then there exist $\eps>0$ and such random vectors 
$X_n\sim\Bex(\pp_n)$ and $X_n'\sim\Bex(\pp_n')$ such that
$
\sumNn\rho(\pni,\pni')< 1/n 
$ 
but $\dtv(X_n,X_n')\ge\eps$, but this contradicts Theorems \ref{T1}(i)
and \ref{T0}.
\end{proof}

\begin{lemma}
  \label{L1ii}
For every constants $C_1,C_2$ and $\eps>0$, there exists $\gd>0$ such
that if $X\sim\Bex(\pp)$ and $X'\sim\Bex(\pp')$ for two deterministic
vectors $\pp=\set{p_i}\xiin$ and $\pp'=\set{p_i'}\xiin$ of the same
length $N\le\infty$, and these satisfy
$
\sumiN\rho(p_i,p_i')\le C_1
$ 
and further, for every $i\in[N]$, 
$p_i\le C_2 p'_i$ and $(1-p_i)\le C_2(1-p'_i)$,
then for every set $A$ with $\P(X'\in A)<\gd$, we have $\P(X\in A)<\eps$.
\end{lemma}

\begin{proof}
  If not, it would be possible to find, for some fixed $C_1$, $C_2$ and
  $\eps$, sequences $X_n$ and $X_n'$ of such random vectors and sets
  $A_n$ such that $\P(X'_n\in A_n)<1/n$ and $\P(X_n\in A_n)\ge \eps$.
In particular, $X_n\not\cont X_n'$.

On the other hand, \eqref{t1iia} and \eqref{t1iib} hold for these
random vectors (since the sums in \eqref{t1iib} vanish for any $C\ge C_2$), 
and thus \refT{T1}(ii) yields $X_n\cont X_n'$, which 
is a contradiction.
\end{proof}

\begin{proof}[Proof of \refT{T2}]
\pfitem{i}
Let $\eps>0$ and choose $\gd>0$ as in \refL{L1i}. Then, by \refL{L1i}
applied to the conditioned variables $\xppn$ and $\xppny$,
if $\sum_i\rho(\pni,\pni')<\gd$, then  $\dtv\bigpar{\xppn,\xppny} <\eps$.
Since 
$$
\P(X_n\in A)-\P(X'_n\in A)=\E\bigpar{
\P(X_n\in A\mid\pp_n)-\P(X'_n\in A\mid\pp_n')}
$$
for every measurable $A\subseteq\cXX_n=\oiset^{N(n)}$, it follows that
\begin{equation*}
  \dtv(X_n,X_n')
\le
\E\dtv\bigpar{\xppn,\xppny}
\le \eps+\P\Bigpar{\sum_i\rho(\pni,\pni')\ge\gd}.
\end{equation*}
The latter probability tends to 0 by assumption, and since $\eps$ is
arbitrary, this yields $\dtv(X_n,X_n')\to0$.

\pfitem{ii}
Let $(A_n)_n$ be an arbitrary sequence measurable sets with
$A_n\subseteq\cXX_n=\oiset^{N(n)}$ and 
let $\eps>0$. 

By \eqref{t2iia}, there exists $C_1$ such that 
$\P\bigpar{\sum_i\rho(\pni,\pni')>C_1}<\eps$ for all $n$.
Similarly, by \eqref{t2iib}, there exist $C_2\ge1$ and $n_0$ such that
for $n\ge n_0$,
\begin{equation*}
\P\lrpar{
\sum_{i:\pni>C_2\pni'}\pni
+\sum_{i:\qni>C_2\qni'}\qni>\eps}
<\eps;
  \end{equation*}
in the sequel we consider only $n\ge n_0$.

Define $\pp_n''=\set{\pni''}\xiinn$ by
\begin{equation*}
\pp_n''\=
  \begin{cases}
	\pp_n', & \sum_i \rho(\pni,\pni')>C_1
\text{ or }
\sum_{i:\pni>C_2\pni'}\pni
+\sum_{i:\qni>C_2\qni'}\qni>\eps;
\\
 \pp_n, & \text{otherwise}.
  \end{cases}
\end{equation*}
By our choices of $C_1$ and $C_2$, $\P(\pp_n''\neq\pp_n)<2\eps$, and
we may thus define $X_n''=(\ini'')\xiinn\sim\Bex(\pp_n'')$ such that 
\begin{equation}
  \label{jesper}
\P(X_n''\neq X_n)\le \P(\pp_n''\neq\pp_n)<2\eps.
\end{equation}

Moreover, by the construction, with $\qni''\=1-\pni''$,
\begin{equation}\label{sofie}
\sum_i \rho(\pni'',\pni')\le C_1
\qquad\text{and}\qquad
\sum_{i:\pni''>C_2\pni'}\pni''
+\sum_{i:\qni''>C_2\qni'}\qni''\le\eps.
\end{equation}

Next, define
$\pni'''=\set{\pni'''}\xiinn$ by
\begin{equation}\label{p'''}
\pni'''\=
  \begin{cases}
	\pni', & \pni''>C_2\pni' \text{ or }
    \qni''>C_2\qni';
\\
 \pni'', & \text{otherwise}.
  \end{cases}
\end{equation}
We can construct $X_n'''\sim\Bex(\pp_n''')$ using  maximal couplings
of $(\ini'''\mid\pp_n''')$ and $(\ini''\mid\pp_n'')$ so that, using
\eqref{p'''} and
\eqref{sofie},
\begin{equation}
\label{ew}
  \begin{split}
\dtv\bigpar{\xppnyyy,\xppnyy}
&\le \sum_i  \dtv\bigpar{(\ini'''\mid\pni'''),(\ini''\mid\pni'')}
\\&
\le \sum_i  
|\pni'''-\pni''|
\\&
\le \sum_{i:\pni''>C_2\pni'} \pni''
+ \sum_{i:\qni''>C_2\qni'} \qni''
\le\eps.	
  \end{split}
\end{equation}
Consequently, 
\begin{equation}
  \label{emma}
\dtv(X_n''', X_n'')\le
\E \dtv\bigpar{\xppnyyy,\xppnyy}
\le \eps.
\end{equation}

Furthermore, by \eqref{p'''},
$\pni'''\le C_2\pni'$ and $\qni'''\=1-\pni'''\le C_2\qni'$
and by \eqref{p'''} and \eqref{sofie},
\begin{equation*}
  \sum_i\rho(\pni''',\pni')
\le
  \sum_i\rho(\pni'',\pni')
\le C_1.
\end{equation*}
We can thus apply \refL{L1ii} to the conditioned variables
$(X_n'''\mid\pp_n''')$ and $(X_n'\mid\pp_n')$ for all values of
$\pp_n'''$ and $\pp_n'$. Consequently there exists $\gd>0$ such that
if $\P(X'_n\in A_n\mid\pp_n')<\gd$,
then
$\P(X'''_n\in A_n\mid\pp_n''')<\eps$.
Hence, using Markov's inequality,
\begin{equation*}
  \begin{split}
\P(X'''_n\in A_n)
&=
\E \P(X'''_n\in A_n\mid\pp_n''')
\le \eps +
\P\bigpar{\P(X'_n\in A_n\mid\pp_n')\ge \gd}
\\&
\le \eps + \gd\qw \E\P(X'_n\in A_n\mid\pp_n')
=
 \eps + \gd\qw \P(X'_n\in A_n).
  \end{split}
\end{equation*}
Using \eqref{jesper} and \eqref{emma}, we thus obtain
\begin{equation*}
  \begin{split}
\P(X_n\in A_n)
&\le
\P(X_n\neq X_n'')
+
\dtv(X_n'', X_n''')
+
\P(X'''_n\in A_n)
\\&
\le
 4\eps + \gd\qw \P(X'_n\in A_n).
  \end{split}
\end{equation*}
If we assume that $\P(X'_n\in A_n)\to0$, it follows that
$\limsup \P(X_n\in A_n)\le4\eps$, and since $\eps$ is arbitrary thus 
$\P(X_n\in A_n)\to0$, which shows that $X_n\cont X_n'$.
\end{proof}

\begin{proof}[Proof of \refC{C1}]
In order to apply \refT{T2}, we reorder $\pijx_{i<j}$ to
$\pnix\xiinn$; we do this
without further comment. 
We also let $\qij\=1-\pij$ and  $\qij'\=1-\pij'$.
  \pfitem{i}
By \eqref{rho11}, \whp{} $\rho(\pij,\pij')\le C_0(\pij-\pij')^2/\pij$
for some $C_0$, and
thus \eqref{c1i} implies \eqref{t2i}, and the conclusion follows by
\refT{T2}(i).
  \pfitem{ii}
Similarly, by \eqref{rho11} again, \eqref{c1ii} implies \eqref{t2iia}.
Moreover, for any $C\ge 2$,
\begin{equation}\label{jb}
\sum_{i:\pij'>C\pij}\pij'
\le
\frac1C \sum_{i:\pij'>C\pij}\frac{(\pij')^2}{\pij}
\le
\frac4C \sum_{i}\frac{(\pij'-\pij)^2}{\pij}
.
\end{equation}
Hence, for any sequence $C_n\to\infty$, \eqref{c1ii} implies that
$
\sum_{i:\pij'>C_n\pij}\pij'\pto0$.
Moreover, for any $C\ge 10$, $C\qij\ge1$ and thus
$\qij'\le C\qij$ for all $i,j$. It follows that
\eqref{t2iiboz} with $\pp$ and $\pp'$ interchanged holds, and thus
\eqref{t2iib} with $\pp$ and $\pp'$ interchanged holds. (The latter is
also easily proved directly using \eqref{jb}.)
Consequently, \refT{T2}(ii) yields $\gnppy\cont\gnpp$. 

\pfitem{iii}
The extra assumptions allow us to interchange $\pp$ and $\pp'$ in the
assumptions. Hence (ii) yields both $\gnpp\contr\gnppy$ and
$\gnppy\contr\gnpp$. 
\end{proof}

\begin{proof}[Proof of \refC{C2}]
  An immediate consequence of \refC{C1}, since now
  $(\pij-\pij')^2/\pij = O(\pij^3)$; note also that the assumption in
  (i) implies $\max_{i,j}\pij=o_p(1)$ and thus $\max_{i,j}\pij<0.9$ \whp.
\end{proof}

\newcommand\AAP{\emph{Adv. Appl. Probab.} }
\newcommand\JAP{\emph{J. Appl. Probab.} }
\newcommand\JAMS{\emph{J. \AMS} }
\newcommand\MAMS{\emph{Memoirs \AMS} }
\newcommand\PAMS{\emph{Proc. \AMS} }
\newcommand\TAMS{\emph{Trans. \AMS} }
\newcommand\AnnMS{\emph{Ann. Math. Statist.} }
\newcommand\AnnPr{\emph{Ann. Probab.} }
\newcommand\CPC{\emph{Combin. Probab. Comput.} }
\newcommand\JMAA{\emph{J. Math. Anal. Appl.} }
\newcommand\RSA{\emph{Random Struct. Alg.} }
\newcommand\ZW{\emph{Z. Wahrsch. Verw. Gebiete} }
\newcommand\DMTCS{\jour{Discr. Math. Theor. Comput. Sci.} }

\newcommand\AMS{Amer. Math. Soc.}
\newcommand\Springer{Springer-Verlag}
\newcommand\Wiley{Wiley}

\newcommand\vol{\textbf}
\newcommand\jour{\emph}
\newcommand\book{\emph}
\newcommand\inbook{\emph}
\def\no#1#2,{\unskip#2, no. #1,} 
\newcommand\toappear{\unskip, to appear}

\newcommand\webcite[1]{
\texttt{\def~{{\tiny$\sim$}}#1}\hfill\hfill}
\newcommand\webcitesvante{\webcite{http://www.math.uu.se/~svante/papers/}}
\newcommand\arxiv[1]{\webcite{arXiv:#1.}}
\newcommand\testtom[1]{%
\def\xxa{#1\relax}\def\xxb{\relax}\ifx\xxa\xxb \else [#1]\fi}

\def\nobibitem#1\par{}


\begin{thebibliography}{99}

\bibitem{BAsc} A.-L.~Barab\'asi \& R.~Albert,
 Emergence of scaling in random networks,
 \emph{Science} \vol{286} (1999), 509--512.

\bibitem[Barbour, Holst and Janson(1992)]{SJI} 
A. Barbour, L. Holst \& S. Janson, 
\emph{Poisson Approximation}. 
Oxford University Press, Oxford, 1992.

\bibitem[Bollob\'as(2001)]{Bollobas}
B. Bollob\'as, 
\book{Random Graphs}, 2nd ed., Cambridge Univ. Press,
Cambridge, 2001.

\bibitem[Bollob\'as, Janson and Riordan(2007)]{SJ178}
B.~Bollob\'as, S. Janson \& O.~Riordan,
The phase transition in inhomogeneous random graphs.
\RSA \vol{31} (2007), 3--122.

\bibitem[Bollob\'as, Janson and Riordan(2007+)]{clustering}
B.~Bollob\'as, S. Janson \& O.~Riordan,
Sparse random graphs with clustering.
In preparation.

\bibitem[Britton, Deijfen and Martin-L\"of(2008+)]{BrittonDML}
T. Britton, M. Deijfen \& A. Martin-L\"of,
Generating simple random graphs with prescribed degree distribution,
\jour{J. Statist. Phys.}, to appear.

\bibitem[Durrett(2003)]{Durrett03}
R.~Durrett,
 Rigorous result for the CHKNS random graph model,
\inbook{Proceedings, Discrete Random Walks 2003 (Paris, 2003)},
eds. C. Banderier \& Chr. Krattenthaler,
{Discrete Mathematics and Theoretical Computer Science}
 \textbf{AC} (2003),
pp. 95--104,
 \webcite{http://dmtcs.loria.fr/proceedings/}

\bibitem[Durrett(2007)]{Durrett}
R. Durrett,
\book{Random Graph Dynamics}.
Cambridge Univ. Press, Cambridge, 2007.

\bibitem[van den Esker, van der Hofstad and Hooghiemstra(2008+)]{EHH}
H. van den Esker, R. van der Hofstad \& G. Hooghiemstra,
Universality for the distance in finite variance random graphs:
Extended version.
\arxiv{math/0605414v2}

\bibitem[Greenhill, Janson, Kim and Wormald(2002)]{SJ140} 
C. Greenhill, S. Janson, J. H. Kim \& N. C. Wormald. 
Permutation pseudographs and contiguity. 
\CPC 
\vol{11}\no3 (2002), 273--298. 

\bibitem{JS}
J.~Jacod \& A. N.~Shiryaev,
\book{Limit Theorems for Stochastic Processes}.
\Springer, Berlin, 1987.

\bibitem[Janson(1995)]{SJ102}
S. Janson,
Random regular graphs: asymptotic distributions and contiguity. 
\CPC 
\vol4 (1995), 369--405.


\bibitem[Janson(2007+)]{SJ195}
S. Janson,
The probability that a random multigraph is simple. 
Preprint, 2006.
\arxiv{math.CO/0609802}


\bibitem[Janson, \L uczak and Ruci\'nski(2000)]{JLR}
S. Janson, T. \L uczak \& A. Ruci\'nski,
\book{Random Graphs},
\Wiley, New York, 2000.

\bibitem{KW} S.~Kalikow \& B.~Weiss,
 When are random graphs connected?
 \emph{Israel J. Math.} \vol{62} (1988), 257--268.

\nobibitem{Kallenberg}
O. Kallenberg,
\book{Foundations of Modern Probability.}
2nd ed., Springer, New York, 2002. 

\bibitem[Oosterhoff and van Zwet(1979)]{OZ}
J. Oosterhoff \& W. R. van Zwet,
A note on contiguity and Hellinger distance. 
\inbook{Contributions to statistics}, pp. 157--166, 
Reidel, Dordrecht, 1979.

\bibitem[Riordan(2005)]{Rsmall}
O.~Riordan,
The small giant component in scale-free random graphs,
\CPC \vol{14} (2005), 897--938.

\bibitem{Wormald}
N. C. Wormald, 
Models of random regular graphs.  
\inbook{Surveys in combinatorics, 1999 (Canterbury)},  239--298, 
London Math. Soc. Lecture Note Ser., 267, 
Cambridge Univ. Press, Cambridge, 1999. 


\end{thebibliography}
\end{document}